\documentclass[10pt]{article}

\usepackage{graphics}
\usepackage{graphicx}

\usepackage[a4paper, left=35mm,right=35mm,top=34mm,bottom=34mm]{geometry}
\usepackage[utf8]{inputenc}
\usepackage[T1]{fontenc}
\usepackage[english]{babel}

\usepackage{enumerate}
\usepackage{graphicx}
\usepackage{hyperref}
\hypersetup{
    colorlinks=true,
    linkcolor=blue,
    filecolor=magenta,      
    urlcolor=cyan,
}
\usepackage{listings}
\usepackage{color}

\definecolor{dkgreen}{rgb}{0,0.6,0}
\definecolor{gray}{rgb}{0.5,0.5,0.5}
\definecolor{mauve}{rgb}{0.58,0,0.82}

\usepackage{mathtools,amsthm,amssymb,amsfonts}
\usepackage{algorithm}
\usepackage{algorithmic}
\usepackage{caption}
\usepackage{subcaption}
\makeatother
\theoremstyle{plain}

\theoremstyle{definition}

\theoremstyle{remark}

\usepackage{caption} 
\usepackage{fancyhdr}
\author{
  {\normalsize Marien Chenaud}\thanks{MICS, CentraleSupélec, Université Paris Saclay, 3 rue Joliot Curie, 91190 Gif-sur-Yvette, France.} \thanks{Transvalor S.A. E-Golf Park, 950 Avenue Roumanille, 06410 Biot, France.} 
  \and
  {\normalsize Fr\'ed\'eric Magoul\`es}\thanks{Corresponding author, frederic.magoules@hotmail.com} \footnotemark[1]
  \and
  {\normalsize José Alves}\footnotemark[2] 
}
\title{Physics-Informed Graph-Mesh Networks for PDEs: A hybrid approach for complex problems}

\date{}

\begin{document}
\maketitle
\thispagestyle{fancy}

\begin{abstract}
The recent rise of deep learning has led to numerous applications, including solving partial differential equations using \textit{Physics-Informed Neural Networks}. This approach has proven highly effective in several academic cases. However, their lack of physical invariances, coupled with other significant weaknesses, such as an inability to handle complex geometries or their lack of generalization capabilities, make them unable to compete with classical numerical solvers in industrial settings. In this work, a limitation regarding the use of automatic differentiation in the context of physics-informed learning is highlighted. A hybrid approach combining physics-informed graph neural networks with numerical kernels from finite elements is introduced. After studying the theoretical properties of our model, we apply it to complex geometries, in two and three dimensions. Our choices are supported by an ablation study, and we evaluate the generalisation capacity of the proposed approach.
\end{abstract}

\section{Introduction}
\label{sec:intro}

Partial differential equations (PDEs) arise in almost every known physical system and
at every scale, whether looking at macro scales with continuous mechanics (fluid or solid), the very small with quantum physics or the very large with astrophysics. Differential equations constitute also the basis of other domains such as finance or population growth in biology. In this context, and since the exact solutions of these equations are generally not reachable, an efficient and
accurate technique to estimate these solutions is critical. In the last century, the theory and numerical procedures tackling these problems have been developed thoroughly. One of the most successful techniques is the \textit{finite element method} (FEM), which has proven to be very efficient numerically. However, depending on the model size and physical non-linearities, such models can incur in considerable computational costs. In industrial settings, with complex geometries and various physics interacting with each other, some efforts
need to be made in order to balance model complexity and computational requirements. Consequently, considering the success of deep learning in addressing a broad range of problems efficiently, it is reasonable to pursue the development of novel deep learning techniques for physical applications. However, since accurate numerical simulations can be expensive in either (or both) run-time and engineering design time, data acquisition can be expensive. Real experimental data can be even more expensive than the simulation data. As a result, available datasets are rare and scarce. Classical deep learning frameworks may thus be inapropriate for this range of problems, therefore, the objective is to develop new techniques combining physical knowledge and the scarce available data to create accurate and trustworthy models.
\newline

Some early attempts, such as \cite{lagaris1998artificial} have been conducted to train deep learning models directly on the physical knowledge we have of the problem, instead of a data-driven approach. This subject has gained attention since the seminal work of \cite{raissi2019physics} and their physics-informed neural networks (PINNs). Since then, many works have been focusing on this task, with many variations and applications \cite{el2023deep,samaniego2020energy,geneva2020modeling,sheng2022pfnn}. However, few of them have been able to use this approach to tackle complex, three-dimensional, industrial problems, and the experiments are often conducted on academic cases, for which traditional numerical techniques are suitable. In this work, after a quick presentation of the PINNs framework and their possible extension to graph networks, along with a justification of the interest of this approach from a physical and numerical perspective, we will point out a limitation with the current way of computing PDE residuals by autodifferentiation. To solve this problem, which limits the ability of PINNs to handle complex geometries, we then propose a new framework, combining a numerical solver to compute physical residuals, and a physically suitable architecture for our Physics-Informed Graph-Mesh Network model. The subsequent investigation is then carried out with a numerical study, first with a two-dimensional mesh to illustrate the limitations of autodifferentiation, and then on a complex, three-dimensional case, to demonstrate the extrapolation ability of our model. This work follows the conference article \cite{chenaud2023physics}. Significant improvements have been made regarding the variational formulation of the problem, the model's architecture and the consideration of physical invariances. Moreover, the generalization properties of our model have been investigated.

\subsection{Partial differential equations and PINNs}
\label{subsec:pinn}

Here, we quickly present the physics-informed framework, as introduced in \cite{raissi2019physics}. This overview has been adapted from \cite{chenaud2023physics}. \newline

For the sake of simplicity, the following section will focus on steady state formalism. A steady-state PDE can be written in a general form like Equation \eqref{eq:steadyPDE}.

\begin{subequations}
    \begin{equation}
        \mathcal{E}(u)(x) = 0 \quad  \forall x \in \Omega, \quad \mathcal{E}\mbox{ being any differential operator,}
    \end{equation}

    \begin{equation}
        \mathcal{B}(u)(x) = 0 \quad \forall x \in \partial \Omega \quad \mbox{(Boundary conditions)}.
        \label{eq:BCsteady}
    \end{equation}
    \label{eq:steadyPDE}
\end{subequations}

In this formulation, $\mathcal{E}$ represents a general, differential operator. $\Omega$ is the domain in which the equation has to be solved (usually, an open, connected set of $\mathbb{R}^d$ for any integer $d \geq 1$). The solution we look for is $u$, function of the variable
space $x \in \Omega$. Finally, $\mathcal{B}$ is any boundary operator (Dirichlet, Neumann or mixed usually,
but it can be more general). In the following, we will suppose that this system of equation is well-posed, meaning that there exists a unique solution to it, which continuously depends on the problem parameters.
\\

The original idea to find an approximation of the solution of steady PDEs such
as \eqref{eq:steadyPDE} with neural networks was to use a rather simple network architecture,
the Multi-Layer Perceptron (MLP).
This model produces a predicted solution $\hat{u}$, and in order to ensure that the PDE is respected in the domain by this solution, some \textit{collocation points} $(x_i)_{i = 1,...,N}$
are used. These are points chosen, randomly or not, in the domain $\Omega$, for which the
residuals $\mathcal{E}(\hat{u})(x_i)$ are computed.
The loss function that needs to be minimized during training is hence:

\begin{equation}
    \mathcal{L}(\hat{u}) = \frac{1}{N} \sum^N_{i=1} ||\mathcal{E}(\hat{u})(x_i)||^2 := \mathcal{L}_r(\hat{u}).
\end{equation}

A neural network with smooth activation functions is differentiable, so   $\mathcal{E}(\hat{u})$ can be computed, either, for simple cases such as in \cite{lagaris1998artificial}, analytically, or with auto differentiation techniques provided by frameworks such as Pytorch \cite{pytorch2019}. The minimization procedure of the loss is typically done by generalized gradient descent techniques such as Adam \cite{kingma2014adam} or L-BFGS. To ensure that the boundary conditions are respected as well, some collocation
points are selected on the boundary too, and the error on the boundary conditions at these points is included in the loss function. The loss is hence:
\begin{equation}
    \mathcal{L}(\hat{u}) = \lambda_{r} \mathcal{L}_r(\hat{u}) + \lambda_{u_b}\mathcal{L}_{u_b}(\hat{u}),
    \label{eq:loss-relaxedBC}
\end{equation}

with $\lambda_{r}, \lambda_{u_b}$ some hyperparameters that need to be tuned. Here, supposing there are $M$ data points
located on the boundary $\partial \Omega$ for which the expected value $u_b$ is known, the loss term $\mathcal{L}_{u_b}$ is:
\begin{equation}
    \mathcal{L}_{u_b}(\hat{u}) = \frac{1}{M} \sum^M_{i=1} ||\hat{u}(x_i) - u_b(x_i)||^2.
    \label{eq:loss-BC}
\end{equation}

There are many different ways to select the critical hyperparameters $\lambda_{r}$ and $\lambda_{u_b}$, see, for instance, \cite{wang2021understanding,berg2018unified,leake2020deep}. \newline

Models trained in this framework are known as \textit{physics-informed neural networks} (PINNs). Many variations of this framework have been studied in the litterature, to overcome various pathologies such as the difficulty to model high frequencies, or competing terms in the loss. For a more complete overview of these pathologies and some efficient solutions, see \cite{wang2021understanding, wang2021eigenvector, wang2022and}.

\subsection{Deep learning on graphs}
\label{subsec:graphDL}

Many physical, social, chemical and biological processes involve interactions between individual entities. These interactions can usually be described using graph-based models, hence the increasing interest of these models in deep learning. For an extensive review and some examples of these methods, see \cite{wu2020comprehensive,zhang2020deep}. These approaches have proven to be very accurate, due to the \textit{inductive biases} of these models. This concept has been studied thoroughly in \cite{battaglia2018relational}. In their work, the authors showed that each machine learning model has inherent biases, which should be taken into account while selecting a model for a given task. The chosen model must align with the properties of the system being studied to ensure its performance. For example, graph-based methods propagate information step by step, from closer to more distant nodes and edges. These models exhibit a locality bias and node-edge permutation invariances which can be correlated with social interconnections between social network users, which explains their use in such applications.
\newline

There are many ways to build a graph network, the key point being how the information diffusion through the graph (nodes and edges) is performed \cite{bianchi2021graph,kipf2016semi,defferrard2016convolutional}.

Recently, \cite{battaglia2018relational} developed a general framework for graph-based approaches, allowing to tackle different tasks with simple variations of the same procedure. A graph is considered to be an aggregation of three different kinds of information, namely features of nodes, edges and the graph. The authors present a general message-passing procedure to diffuse the features within the graph through these three entities. In our case, we focus on a node prediction task, with both nodes and edges features (and no graph features). The graph model, called graph network, is made of several \textit{graph network blocks}, each block corresponding to a generalization of a convolutional layer. First, edge features are updated using the features of adjacent nodes and the old edge features. Then, nodes are updated, with the new edge features. Each update is handled by a neural network, usually a multi-layer perceptron (MLP). For the sake of completeness, the main steps of the message-passing operation through these graph network blocks are presented here. For a more extensive analysis of this mechanism, see \cite{battaglia2018relational}. \newline

First, we call a `graph' an object $G = (V,E)$. $V = \{v_i\}_{i=1:N^v}$ is the set of nodes, of cardinality $N^v$, and each $v_i$ is a node's attribute (such as the value of a physical field on this particular node). $E = {(e_k, r_k, s_k)}_{k=1:N^e}$ is the set of edges (of cardinality $N^e$), where each $e_k$ is the edge’s attribute, $r_k$ is the index of the receiver node, and $s_k$ is the index of the sender node. In the applications below, the connectivity of the graph is extracted from the connectivity of the corresponding mesh, created to discretize an object, in the same fashion as what is done in Finite Element Method.\newline

A graph network block is made of two `update' functions $\phi^e$ and $\phi^e$, operating respectively on the edges  and the nodes, and an `aggregation' function $\rho$. In our model, the aggregation function is chosen to be the sum function. The $\phi^e$ is mapped across all edges to compute per-edge updates, and the $\phi^v$ is mapped across all nodes to compute per-node updates. The function $\rho$ take a set as input and reduces it to a single element by summing them, the sum representing the aggregated information. Mathematically, one pass through a graph network block can be modelized by the operations \eqref{eq:gn1} and \eqref{eq:gn2}.

\begin{align}
    \label{eq:gn1}
    e'_k & = \phi^e(e_k, v_{r_k}, v_{s_k}) \quad & \bar{e}'_i & = \rho(E'_i) \\
    v'_i & = \phi^v(\bar{e}'_i,v_i)              &            &
    \label{eq:gn2}
\end{align}

\noindent where $E'_i = \{e'_k, r_k, s_k\}_{r_k=i, k = 1:N^e}$.  \newline

\noindent The updated graph is then $G' = (V', E'),$ with $V' = \{v'_i\}_{i=1:N^v}, \, E' = \{e'_k\}_{k=1:N^e}$. While the aggregation function $\rho$ is fixed, the update functions $\phi$ can be any function, and are therefore chosen to be modelized by a neural network, a MLP, with trainable parameters. During training, these functions will evolve, guiding the diffusion of information through the graph.

\begin{figure}[!h]
    \centering
    \includegraphics[width = \textwidth]{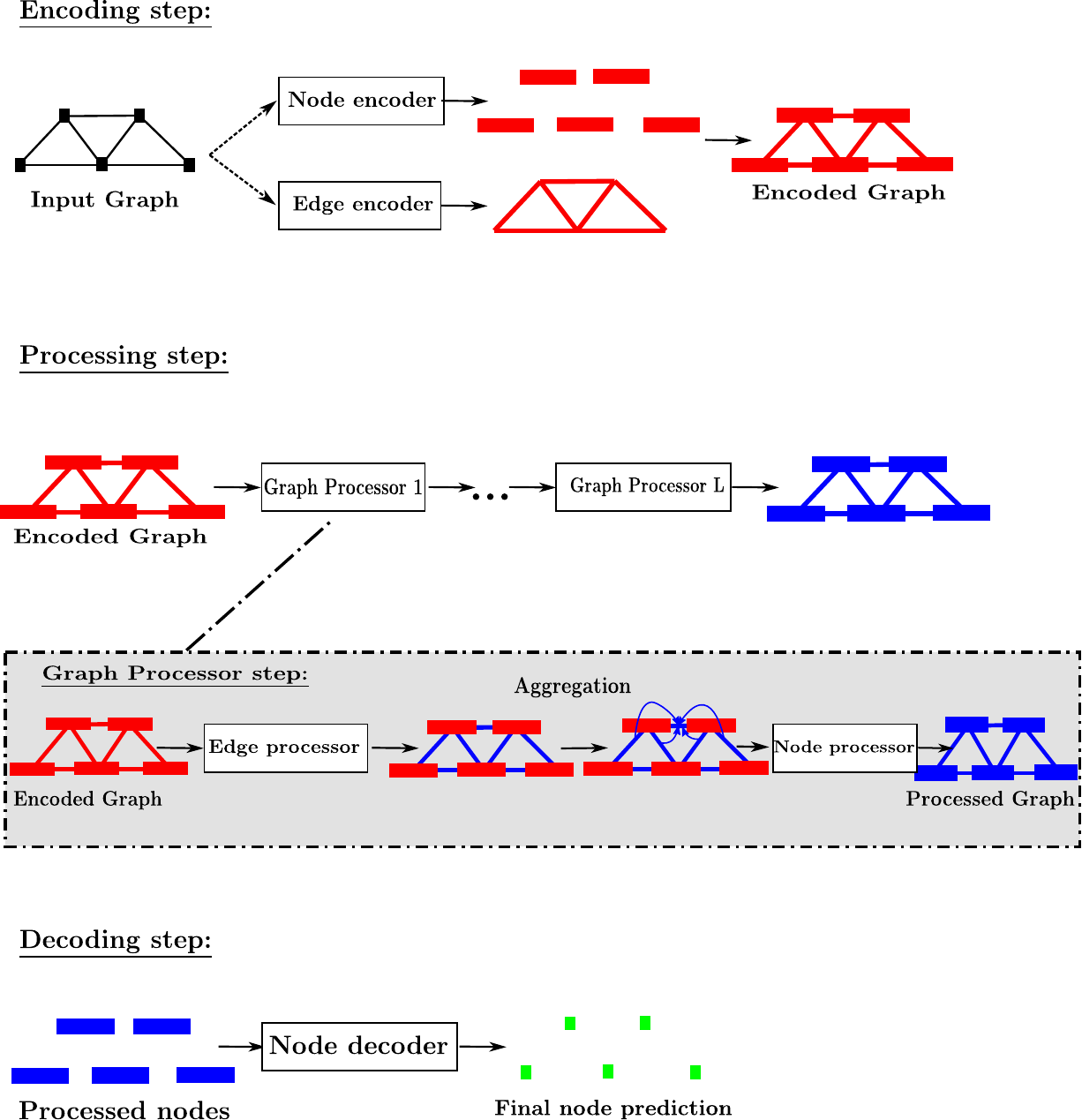}

    \caption{Model architecture. First, an encoder processes the physical input data. Next, this enriched input is fed through $L$ graph network blocks. Finally, a node decoder is applied to predict the target physical field. The black elements refer to the physical input, the red elements are the encoded data, the blue objects are the processed data, and the green nodes represent the output field, i.e the model's prediction.}
    \label{fig:architecture}
\end{figure}

A recent trend involves the use of the previously mentioned graph network, in combination with an encoder-decoder structure. This approach is employed in \cite{pfaff2020learning}, \cite{lam2022graphcast} and \cite{sanchez2020learning}. Both approaches feature a similar model architecture, which can be described as follows. First, the input physical data is fed into an encoder (a MLP). The corresponding output is an embedding of the input in a higher-dimensional space. While the input node fields typically have a dimension of less than 10 (one dimension for each input scalar physical field, 2 or 3 dimensions for each input vector field), the enriched embedding has typically a dimension of 64. Next, the enriched, higher-dimensional encoded data is transmitted through several graph net blocks, or graph processors, which play the role of the message-passing step through the graph. Finally, a decoder is used to attain the desired physical target. This approach seamlessly integrates the inductive biases inherent in graph networks with the expressiveness of neural networks, hence its performance. A similar architecture is used in this work. Figure \ref{fig:architecture} illustrates the key steps of this model architecture.

\subsection{Related works}

In \cite{chamberlain2021grand}, the message-passing step over a graph is explained as a diffusion process, and the authors prove that the different convolution schemes can be seen as different techniques to solve a discretized heat diffusion PDE. This strong link between graph convolutional architectures and physical diffusion processes, along with the natural inductive biases explained in \cite{battaglia2018relational}, motivates the use of this framework to simulate physical processes. The capacity of graph networks to simulate complex, dynamic physical problems has also been proven in \cite{sanchez2020learning}, where different fluid-solid interactions are simulated using a graph network. As a result, several recent works have been leveraging this approach to solve PDEs.
\newline

In \cite{gao2022physics}, the authors take advantage of the natural link between graphs and meshes to tackle the problem of solving static PDEs. A mesh is seen as a graph, and the neural network learns the discretized solution of the variational problem in each point of the graph. However, only two-dimensional cases are tackled. \cite{belbute2020combining} incorporated a numerical solver inside a graph model. The result of the solver, which was computed on a coarse mesh, is fed into the graph network to help convergence. \newline

\cite{meethal2023finite} used a loss constructed on finite element residuals to train their model, by using a finite element code to compute the stiffness matrix and combining it with their model's prediction. This alleviates the issue of computing residuals by autodifferentiation. The model used is simply a neural network, their approach is not graph-based. Therefore, the generalization capability of the model is not explored. \newline

Finally, to alleviate the issue of computing PDE residuals by autodifferentiation, \cite{xiang2022rbf} approximated the spatial derivatives with finite differences. However, the applications are limited to two-dimensional domains, with no investigation on the generalization capabilities.

\section{Extension of PINNs on complex settings}
\label{sec:customGrad}

\subsection{Physical invariances}
\subsubsection*{Learning physical operators}

Although PINNs have demonstrated the capability to solve different direct and inverse problems \cite{raissi2019physics,samaniego2020energy}, they typically require a long training time in comparison to conventional numerical methods such as finite element methods, and they lack generalization abilities. Therefore, a slight change in the limit conditions or the studied domain requires re-training the model. To adress this behaviour, several works have tried to learn physical operators instead of direct solutions, and then to obtain the solutions by applying these physical operators to a particular setting \cite{wang2021learning,podina2023universal,muller2023exact,greydanus2019hamiltonian}. With this approach, the basic physical properties one would expect from a system, such as conservation laws, invariances or equivariances are inherently respected by the model's prediction. Another way to guarantee this behaviour is by using a suitable model's architecture. \newline

\subsubsection*{Physical invariances and model architecture}

\begin{figure}[h]
    \centering
    \begin{subfigure}[b]{0.35\textwidth}
        \centering
        \includegraphics[width=\textwidth]{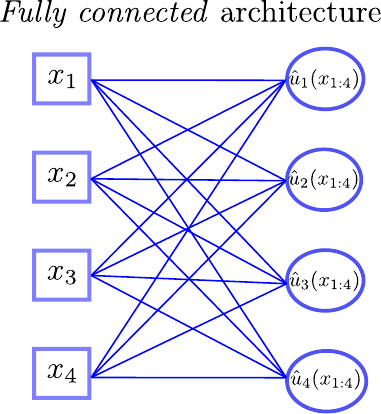}
        \caption{Fully connected architecture.}
        \label{fig:fullyconnected}
    \end{subfigure}

    \vspace{1cm}

    \begin{subfigure}[b]{0.9\textwidth}
        \centering
        \includegraphics[width=\textwidth]{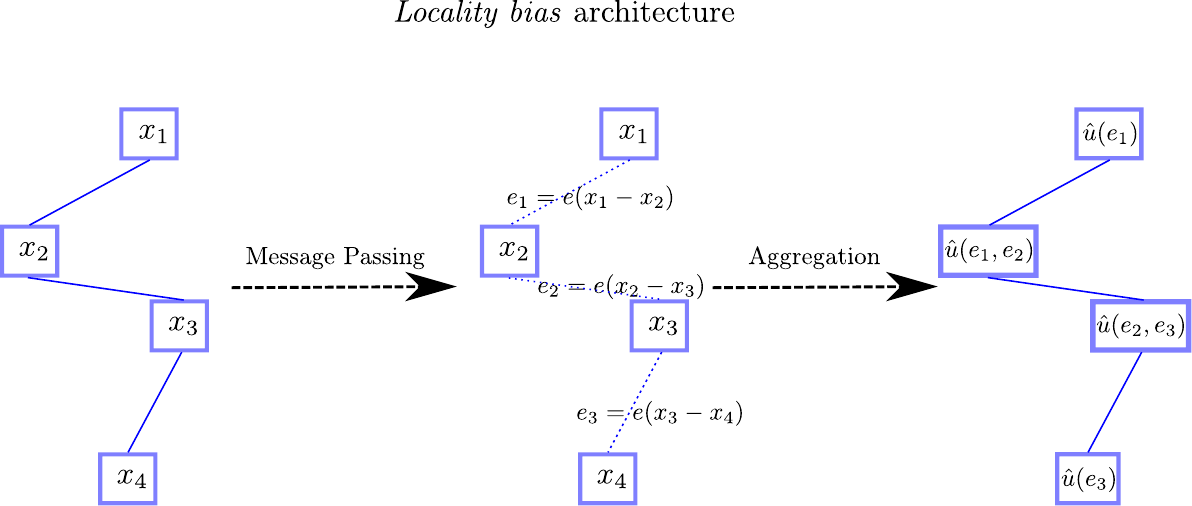}
        \caption{Locally biased architecture.}
        \label{fig:localbias}
    \end{subfigure}

    \caption{Fully connected vs locally biased architectures. The output formula show a translation invariance for the locally biased architecture, as opposed to the fully connected model.}
    \label{fig:connected_bias}
\end{figure}

To be accurate and have good generalization capabilities, the model studied needs to respect the physical invariances of the underlying problem. The model's architecture should account for these invariances, or equivariances. For example, if the problem considered is invariant to space translations, which is the case for most physical problems (where the space origin is arbitrarily set), then the domain coordinates are not suitable as input for the surrogate model. Instead, for a mesh-based approach, relative positions between nodes of the considered geometry are a more suitable input. Figure \ref{fig:connected_bias}, adapted from the discussions in \cite{battaglia2018relational}, illustrates an example of the locality inductive bias. In this general representation, the inputs are denoted $x$, and the output $\hat{u}$. In the case of the locally biased architecture, the edges have attributes given by a function $e$. This function is fairly general, the only restriction being that for an edge between two nodes $x_i$ and $x_j$, the edge attribute should be a function of $x_i - x_j$. For instance, by denoting $d_{ij}$ the Euclidian distance between nodes $x_i$ and $x_j$ and $\mathbin\Vert$ the concatenation operator, $e_{ij} = (x_i - x_j)\mathbin\Vert d_{ij}$ is a suitable edge attribute. While, for graph-based approaches with a locality bias, a perturbation of the input positions in the form $x \mapsto x + \delta$ for any quantity $\delta$ would not lead to a perturbed model output, this is not the case for the fully connected architecture, which is not translation invariant in space.

Most PINN-based approaches disregard this limitation, since input coordinates are typically required to compute the PDE residuals by autodifferentiation. In this work, we propose an approach that respects this physical invariance, which allows to train our model on a given, three-dimensional geometry, and to evaluate it on a very different one, which, to the best of our knowledge, has not been done with PINN-based approaches.

\subsection{A limitation with residuals computations in a physics-informed framework}
\label{subsec:pbGrad}

As explained in section \ref{subsec:pinn}, the physics-informed framework relies on the computation of PDE residuals by the autodifferentiation techniques \cite{pytorch2019} usually used in deep learning for the gradient-based training procedures. This powerfull tool allows to trace back any operation made inside a deep learning model, and to compute its derivative with the chain rule. Therefore, if a position is given as input to the model, applying this tool will allow to compute physical gradients, and then the PDE residuals. However, there are two main issues with this framework.

\begin{itemize}
    \item As seen previously, the spatial coordinates do not respect the invariances required by general conservation principles. Training a model with these coordinates as inputs will prevent its generalization to a slightly different problem, for instance if the geometry is translated in space, or if the geometry is morphed.

    \item The computation of physical derivatives with this framework can be biased. We prove this issue on a simple case for a linear model, since even in this setting the aforementionned bias can arise.
\end{itemize}

Suppose the existence of a linear model $M$, which takes as input a field $\varphi \in \mathbb{R}$ and the position $x = (x_1, \dots, x_n) \in \mathbb{R}^n$ on a given geometry $\Omega \subset \mathbb{R}^n$, and which outputs a scalar field $M(\varphi, x)$. Let $(b_1, \dots, b_n)$ be the canonical basis of $\mathbb{R}^n$. Suppose that the input field $\varphi$, which can be any physical information on the problem, such as the boundary conditions or a physical field at a previous time-step, depends on the position: $\varphi = \varphi(x)$.
\newline
Since $M$ is linear, there exist scalar numbers $\alpha, \gamma$ and a vector $\beta = (\beta_1, \dots, \beta_n) \in \mathbb{R}^n$ such that $M(\varphi, x) = \alpha \varphi(x) + \beta^T \cdot x + \gamma.$ Thus, the corresponding derivative of the output with respect to $x_k$, the $k$-th component of the variable $x$, is:

\begin{equation}
    \nabla M \cdot b_k := \frac{\partial M}{\partial x_k} = \alpha \frac{\partial \varphi}{\partial x_k} + \beta_k.
\end{equation}

In general, physical fields on complex geometries do not have an analytical expression depending on $x$. Hence, in this case where $\varphi$ has not been constructed analytically, there will be no recorded operation in the computational graph between $\varphi$ and $x$. This issue is not limited to any autodifferentiation framework such as Pytorch, and arises whenever derivatives are constructed by tracing back the operations made to construct the final tensors \cite{baydin2018automatic}. The derivative of the output $M(\varphi, x)$ with respect to $x_k$ will therefore be computed once again as $\frac{\partial M}{\partial x_k} \, \underset{\mbox{autodiff}}{=} \, \alpha \frac{\partial \varphi}{\partial x_k} + \beta_k.$ However, $\varphi$ being a direct input to the model, there is no recorded link in the computational graph between $\varphi$ and $x$: therefore, $\frac{\partial \varphi}{\partial x_k} \underset{\mbox{autodiff}}{=} 0 $. Hence, the derivative of $M$ with respect to $x_k$ will be computed as:
\begin{equation}
    \frac{\partial M}{\partial x_k} \, \underset{\mbox{autodiff}}{=} \, \beta_k.
\end{equation}

This issue prevents the use of the physics-informed framework for many general settings, and it would explain why there is a clear lack of work in this direction.

\subsection{Residuals computation with differentiable numerical kernels}
\label{subsec:computeGrad}

As proven in section \ref{subsec:pbGrad}, autodifferentiation may fail to compute physical derivatives in given settings, for instance when additional inputs are given to the model. However, several numerical methods, including finite element methods, have the ability to evaluate these quantities. Finite difference methods could also compute them, but they are not suited for complex and non-structured meshes. Still, it would be troublesome to implement finite element methods within a deep learning framework such as Pytorch. One could call another method outside of the deep learning process during training, but the resulting quantities would be untraceable by the autodifferentiation framework, making it useless for the training procedure. \newline

\begin{figure}[!h]
    \centering
    \includegraphics[width = \textwidth]{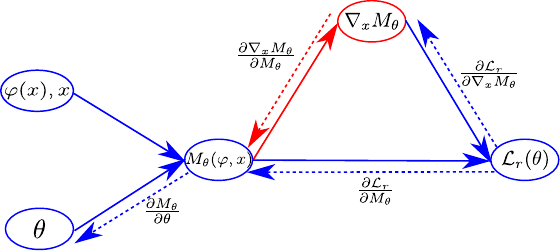}

    \caption{Workflow of the proposed process. The red quantities are computed by the numerical gradient kernel, while the blue quantities are computed inside the deep learning framework. Plain arrows trace the operations made during the forward call, while dashed arrows represent the operations made during the backward pass.}
    \label{fig:workflow}
\end{figure}

Here, we propose to include the call to an external numerical procedure into the training loop, by making this call differentiable. The missing step to do so is the computation of the derivative of the numerical gradient operator with respect to the position. As a matter of fact, since physics-informed models incorporate PDE residuals inside the loss functions, the computation of these residuals, and therefore the computation of physical operators such as gradients or integral computations should be traceable by the autodifferentiation method. Therefore, if one decides to call an outside numerical operator, the derivative of this call should be computed too, to perform backpropagation during training. If, given a physical field predicted by any deep learning model $M_{\theta}(\varphi, x)$, the numerical operator computes the quantity $\nabla_x M_{\theta}(\varphi, x)$, we should be able to compute the following quantity: $\frac{\partial \nabla_x M_{\theta}(\varphi, x)}{\partial  M_{\theta}(\varphi, x)}$. However, numerical operations inside a finite element gradient kernel can be reduced to matrix-vector multiplications or other algebraic operations, and therefore are analytically differentiable. In our proposed approach, these operations are recorded and differentiated, in order to provide the desired derivatives in a custom backward operation inside the autodifferentiation computational graph, allowing to perform gradient-based optimization during training. This procedure is summarized in figure \ref{fig:workflow}.

It is worth noting that this procedure is not specific to neural networks or graph models and could be generalized directly to any machine learning model, denoted $M_{\theta}$ in figure \ref{fig:workflow}. Moreover, this approach is not limited to gradient computations and could be generalized to any physical operator, as long as it is differentiable and that its derivative can be computed numerically. In the following, integrals computations will also be handled by this approach.

\section{Proposed model: PiGMeN, Physics-informed Graph-Mesh Network}
\label{sec:model}
\subsection{Architecture of the model}

In order to benefit from the desired inductive biases of graph networks, along with the expressivity of neural networks, our model architecture follows the ones used in \cite{pfaff2020learning} and \cite{sanchez2020learning}, and described in section \ref{subsec:graphDL}. There are five main components of this model: node and edge encoders, node and edge processors, and node decoders. For each of them, we used multi-layer perceptrons as models. Table \ref{table:architecture} summarizes the architectures of the neural networks for each component. Moreover, a standardization layer is applied before the encoder step, and there are three graph processors (a combination of node and edge processors) in the model.

\begin{table}[ht]
    \centering
    \begin{tabular}{|c||c|c|c|}
        \hline
        Component      & Residual connection & Number of layers & Layer width \\
        \hline
        Node encoder   & No                  & 2                & 64          \\
        \hline
        Edge encoder   & No                  & 2                & 64          \\
        \hline
        Node processor & Yes                 & 2                & 64          \\
        \hline
        Edge processor & Yes                 & 2                & 64          \\
        \hline
        Node decoder   & No                  & 2                & 64          \\
        \hline
    \end{tabular}
    \caption{Hyperparameters of each neural networks used in the global model. For every neural network, the ReLU activation function is used.}
    \label{table:architecture}
\end{table}

The use of residual connections for processors follows the articles \cite{pfaff2020learning,sanchez2020learning}. Our experiments also showed a major improvement when this residual connection is added. Regarding the other parameters, the number of layers and their width, we tried to select the smallest architecture possible, without impacting the performance of the model. We found that increasing the model's complexity, in terms of number of layers or width, did not lead to major improvements in the performances. Therefore, the values of these hyperparameters have been chosen as reported in Table \ref{table:architecture}.

\subsection{Input data for physical invariances}
\label{subsec:input}

As explained previously, the direct spatial coordinates are not a suitable input for a physics-informed model which should be capable of generalization. Therefore, only the relative, spatial coordinates will be used as edge features. For an edge $e_{ij}$ between two nodes $i$ and $j$, with spatial coordinates $x_i$ and $x_j$, the edge coordinate will be the four-dimensional vector $[x_i - x_j \mathbin \Vert d_{ij}]$, where $\mathbin \Vert$ is the concatenation symbol, and $d_{ij}$ the Euclidian distance between nodes $i$ and $j$: $d_{ij} = ||x_i - x_j||$. \newline

As the node features, one could simply provide the boundary conditions field, which would be the null field in $\Bar{\Omega}$, except for non-homogeneous Dirichlet boundary conditions. However, we decided to use instead an approximation of the true solution, which could represent an initial guess of the target field, for instance computed with a coarse numerical solver. More details on this input field, along with a visual representation, are provided in section \ref{subsec:pb}.

\subsection{Training procedure: hybrid computation of loss residuals}
\label{subsec:loss}

Once provided with a reliable method to compute physical derivatives, there are many ways to compute PDE residuals. While \cite{chenaud2023physics} used a direct, strong loss, many works, such as \cite{samaniego2020energy}, show that weak losses, corresponding to energy minimizations, are closely related to the finite element formalism, and therefore are more interpretable physically. Weak formulations allow to incorporate Neumann boundary conditions naturally, and they necessitate lesser-order derivatives computations compared to the direct approach. \cite{yu2018deep} use weak losses to address high-dimensional problems, with conclusive results for problems in dimensions up to 10. Such problems seem unreachble for PINNs with a strong loss formulation. However, including constraints such as boundary conditions in this framework could be challenging, specially for complex geometries. While \cite{samaniego2020energy} strongly enforced boundary conditions with an explicit, algebraic formulation, this would not be possible in our case. Therefore, we propose to adapt a classical technique used in constrained optimization: the Lagrangian formulation. \newline

Since the energy formulation made in \cite{samaniego2020energy} can be seen as a minimization problem, and that the non-zero Dirichlet boundary condition is an additional constraint to the optimization problem, the Lagrangian formulation of this optimization problem, given an energy functional $\mathcal{E}(u)$, would therefore be:

\begin{equation}
    \mathcal{L}(u, \lambda) = \mathcal{E}(u)
    + \int_{\partial \Omega} \lambda (u - u_{\partial \Omega}) dS.
\end{equation}

In this formulation, $u_{\partial \Omega}$ represents the target boundary field, and $\lambda$ is the Lagrange multiplier, associated to the constraint on the minimization problem. The first term accounts for the energy formulation, while the second term accounts for the additional, boundary condition constraint. For instance, the energy formulation derived from the electrostatic PDE \eqref{eq:staticPb}, presented in section \ref{subsec:pb}, would be:

\begin{equation}
    \mathcal{E}(u) = \int_{\Omega}||\nabla u || ^2 \, d\Omega.
\end{equation}

We know, from classical optimization theory, that a necessary condition for a solution $u$ of the constrained minimization problem is that there exists a Lagrange multiplier $\lambda$ such that $(u, \lambda)$ is a saddle point of $\mathcal{L}$. There are ways to directly solve the dual optimization problem with the Lagrangian formulation within a deep learning framework, see, for example, \cite{fioretto2021lagrangian}. However, since in our case we combine this formulation with a data-driven approach, we have only extracted the necessary condition. In order to benefit both from limited available data and the physical knowledge of the problem, the norm of the Lagrangian derivative seems to be a good physical quantity to guide the training. From now on, our models will therefore predict two scalar fields depending on the model's parameters $\theta$, namely $u_{\theta}$ and $\lambda_{\theta}$. The associated loss (or cost) function will then be:

\begin{equation}
    \mathcal{C}(\theta) = \left | \left | \frac{\partial \mathcal{L}(u_{\theta}, \lambda_{\theta})}{\partial u} \right | \right |^2 + \left | \left |\frac{\partial \mathcal{L}(u_{\theta}, \lambda_{\theta})}{\partial \lambda} \right | \right |^2 + \mathcal{C}_{data}(u_{\theta}),
\end{equation}
where $\mathcal{C}_{data}(u_{\theta})$ is the mean absolute error (MAE) between the neural network's prediction $u_{\theta}$ and the true solution $u$:

\begin{equation}
    \mathcal{C}_{data}(u_{\theta}) = ||u - u_{\theta}||_{L^1}.
\end{equation}

This hybrid approach has the benefit of being generalizable to many settings, while usually taking the boundary conditions into account inside the loss can be troublesome. A hyperparameter balancing the physics-informed terms and the data term $\mathcal{C}_{data}$ could be computed, with approaches similar to the one described in \cite{wang2021understanding}, but our numerical experiments have shown that keeping the loss term as is leads to satisfying results, that are not further improved by the tuning of this hyperparameter. Therefore, the numerical experiments have been conducted with this more direct approach. In the future, a method to include all the data knowledge we have on a problem inside the Lagrangian formulation should be investigated. Finally, the integrals and the physical derivatives can be computed accurately thanks to the framework described in section \ref{subsec:computeGrad}.

\section{Numerical experiments in 2D}
\subsection{Presentation of the geometry}
\label{subsec:geometry2d}
\begin{figure}[ht]
    \centerline{
        \begin{subfigure}[b]{0.55\linewidth}
            \includegraphics[width=\linewidth]{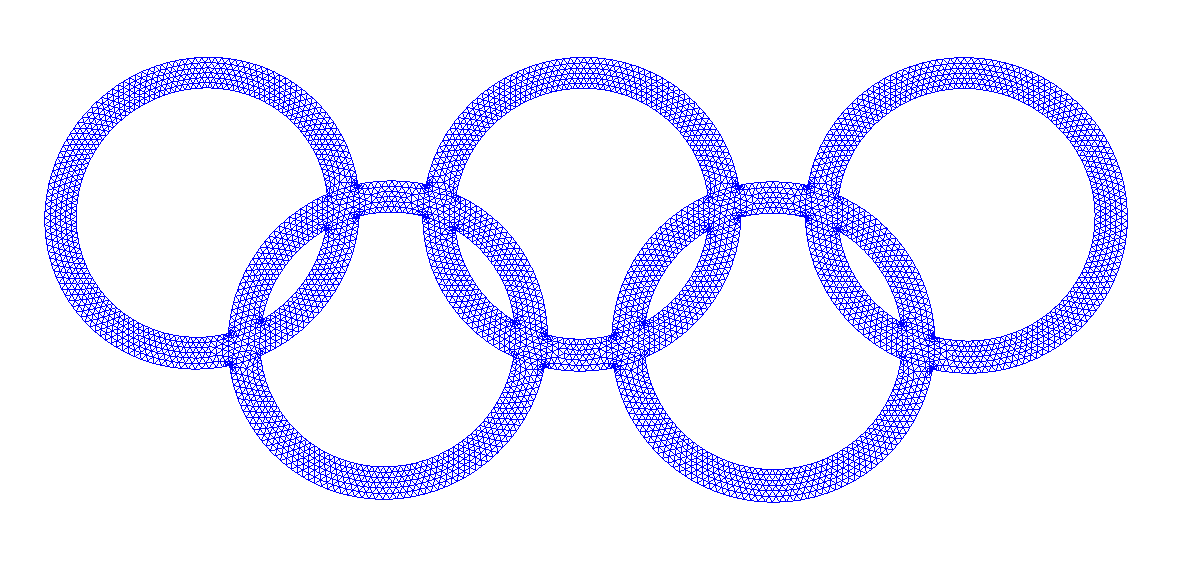}
            \caption{Fine mesh}
        \end{subfigure}

        \begin{subfigure}[b]{0.55\linewidth}
            \includegraphics[width=\linewidth]{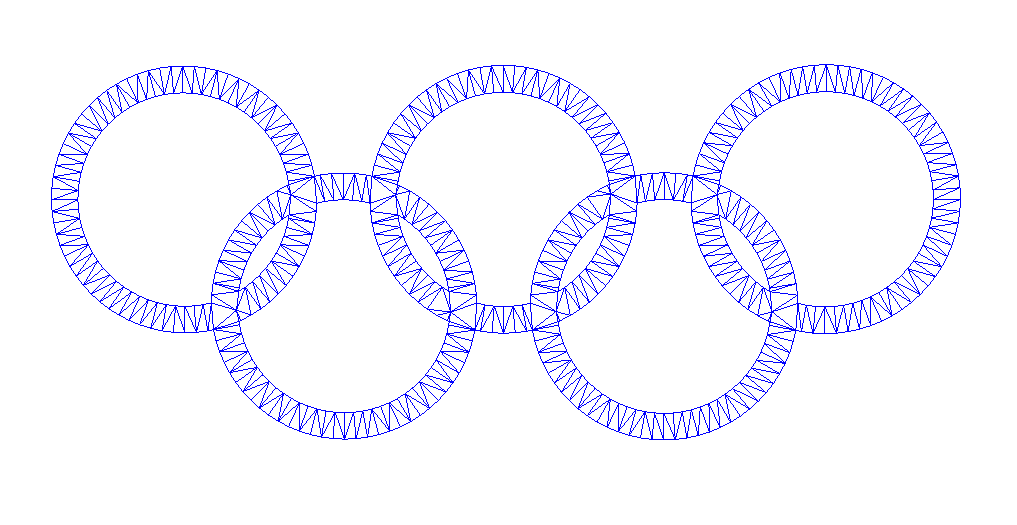}
            \caption{Coarse mesh}
        \end{subfigure}
    }
    \centerline{
        \begin{subfigure}[b]{0.55\linewidth}
            \includegraphics[width=\linewidth]{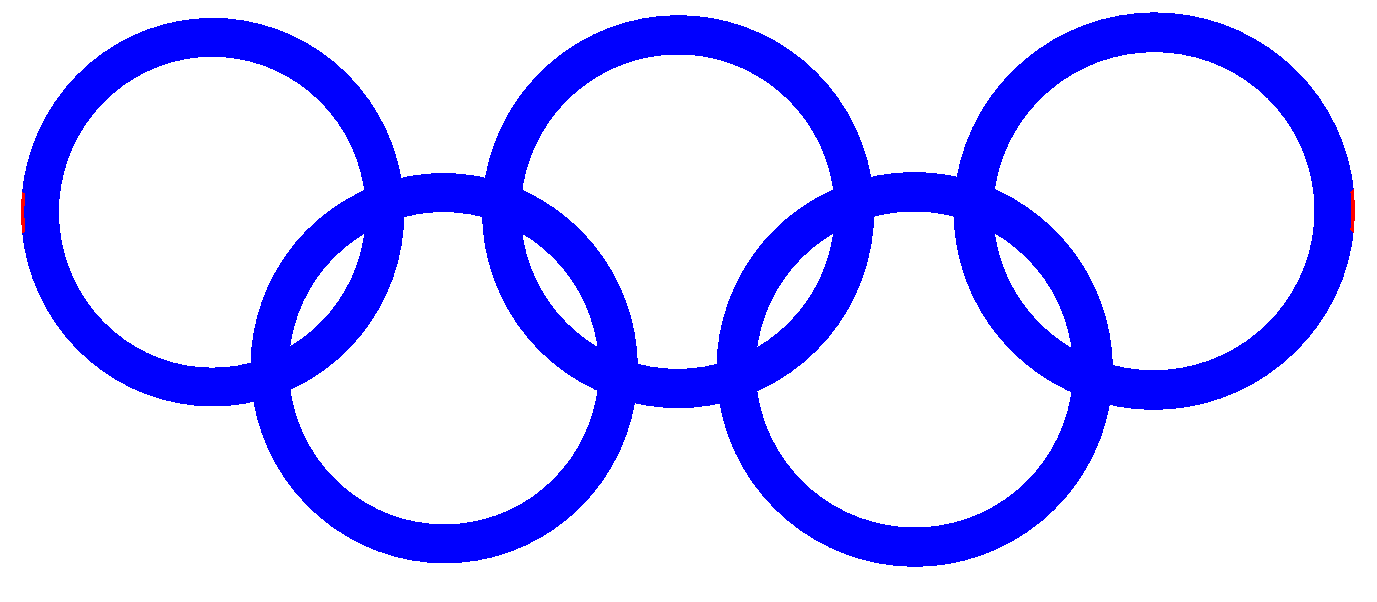}
            \caption{Boundary condition nodes}
            \label{subfig:bc}
        \end{subfigure}
    }
    \caption{Geometry and meshes used for the 2D numerical experiments. The boundary condition nodes are in red.}
    \label{fig:geometry2D}
\end{figure}

To investigate the performance of our approach, a first two-dimensional case has been chosen on a domain $\Omega$. This irregular geometry, chosen to represent the Olympics rings, is challenging, because of its strong non-convexity. Two meshes have been used to investigate the accuracy of the model with respect to the number of nodes. The fine mesh is composed of 5104 nodes and 26520 edges, forming 8840 elements. The coarse mesh is made of 664 nodes with 2064 edges, forming 688 elements. Figure \ref{fig:geometry2D} illustrates the geometry and the chosen meshes.

\subsection{An example of the failure of autodifferentiation}

To illustrate the theoretical result of Section \ref{subsec:pbGrad}, a first case, of a diffusion problem, on the geometry described previously, has been conducted. The considered problem is an electrostatic problem:

The unknown is the electric potential $V$ on a domain $\Omega$. Two portions of the boundary, $\Gamma_1$ and $\Gamma_2$, are fixed, as boundary conditions, and the PDE that needs to be solved for $V$ is:

\begin{subequations}
    \begin{equation}
        \Delta V(x) = 0 \quad \forall \, x \in \Omega,
    \end{equation}

    \begin{equation}
        V(x) = 0 \quad \forall \, x \in \Gamma_1,
    \end{equation}

    \begin{equation}
        V(x) = 100 \quad \forall \, x \in \Gamma_2,
    \end{equation}
    \begin{equation}
        \frac{\partial V}{\partial n}(x) = 0 \quad \forall \, x \in \partial \Omega \backslash (\Gamma_1 \cup \Gamma_2).
    \end{equation}
    \label{eq:staticPb2d}
\end{subequations}

The subset $\Gamma_1$ is located on the right side of the domain, and $\Gamma_2$ on the left side (see Figure \ref{subfig:bc}).

\begin{figure}[ht]
    \centering
    \includegraphics[width = \linewidth]{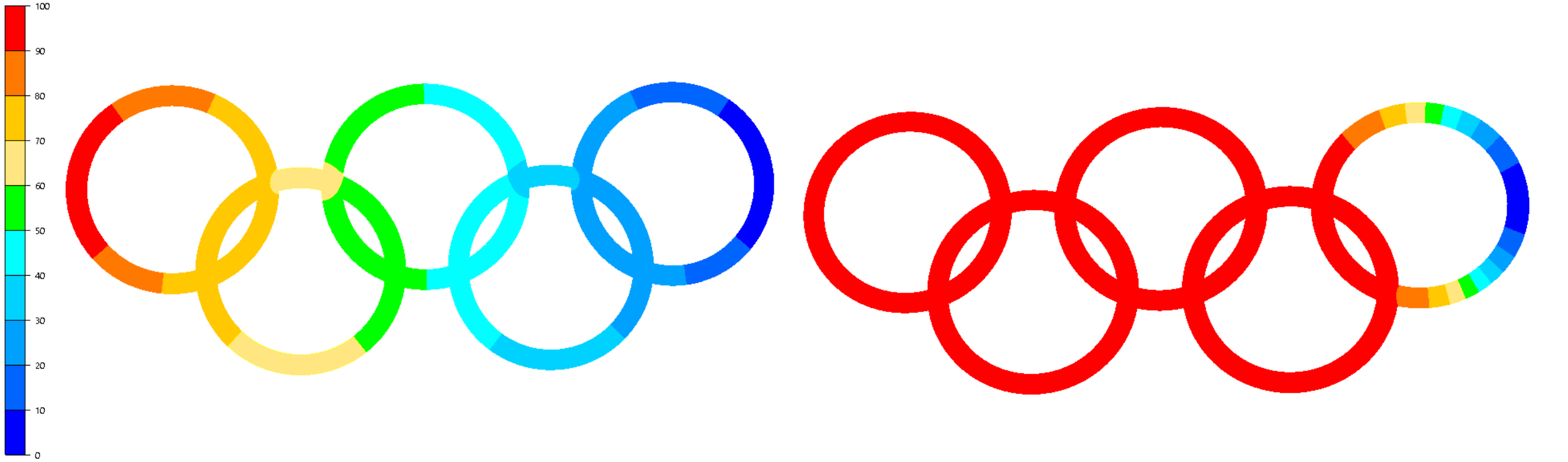}

    \caption{(Left) Target potential, given as input to the model. (Right) Predicted potential.}
    \label{fig:pot2D}
\end{figure}

To make sure that the observed behaviour is caused by the result of Section \ref{subsec:pbGrad}, a simple MLP with two hidden layers and 100 neurons per layer has been used, with a strong formulation of the PDE. The problem has been designed to be particularly simple: as additional input field $\varphi$ to the model, the target solution $V^*$, computed by Finite Element Method, has been added. Therefore, the problem is trivial, since the target solution is given as input to the model. However, while performing the training, this additional input field prevents the accurate computation of the PDE residuals, as pointed out in Section \ref{subsec:pbGrad}. The training has been performed for 5000 epochs, with the Adam optimizer and a learning rate of $5 \times 10^{-3}$. Note that increasing the number of epochs, or changing the training's hyperparameters of MLP structure yielded similar results. Figure \ref{fig:pot2D} plots the target potential, given as input to the model, and the model's prediction. After training, the computed PDE residuals are of the magnitude $2 \times 10^{-5}$.

This simple example illustrates how additional input fields to the model leads to wrong PDE residuals computations, preventing the use of enriched inputs.

\subsection{An elasticity problem}

\subsubsection*{Case presentation}

In this section, our physics-informed model is evaluated on an elastic problem, on the geometry presented in Section \ref{subsec:geometry2d}. The considered PDE, where the unkown $\mathbf{u}$ is the displacement, is as follows:

\begin{subequations}
    \begin{equation}
        \nabla \sigma(\mathbf{u})(x) = 0 \quad \forall \, x \in \Omega,
    \end{equation}

    \begin{equation}
        \sigma(\mathbf{u}) = \lambda tr(\varepsilon) I + 2 \mu \varepsilon, \quad \varepsilon = \frac{1}{2}(\nabla \mathbf{u} + \nabla \mathbf{u}^T),
    \end{equation}

    \begin{equation}
        \mathbf{u}(x) = (1,0) \quad \forall \, x \in \Gamma_1,
    \end{equation}

    \begin{equation}
        \mathbf{u}(x) = (-1,0) \quad \forall \, x \in \Gamma_2,
    \end{equation}
    \begin{equation}
        \frac{\partial \mathbf{u}}{\partial n}(x) = 0 \quad \forall \, x \in \partial \Omega \backslash (\Gamma_1 \cup \Gamma_2).
    \end{equation}
    \label{eq:elasticPb2d}
\end{subequations}

The physical parameters $\lambda$ and $\mu$ are chosen equal to 1. The problem \eqref{eq:elasticPb2d} has been solved by Finite Element method using the meshes presented in Figure \ref{fig:geometry2D}. As explained in section \ref{subsec:geometry2d}, in order to investigate the interpolation ability of the model, two trainings have been performed: one on the fine mesh of Figure \ref{fig:geometry2D}, with 5104 nodes, and another on a much coarser mesh of the same geometry, with 664 nodes.

\subsubsection*{Results and discussion}

\begin{figure}
    \centerline{

        \rotatebox{90}{\hspace{1.2cm} $u_x$}

        \begin{subfigure}[b]{0.55\linewidth}
            \caption{Target displacement}
            \includegraphics[width=\linewidth]{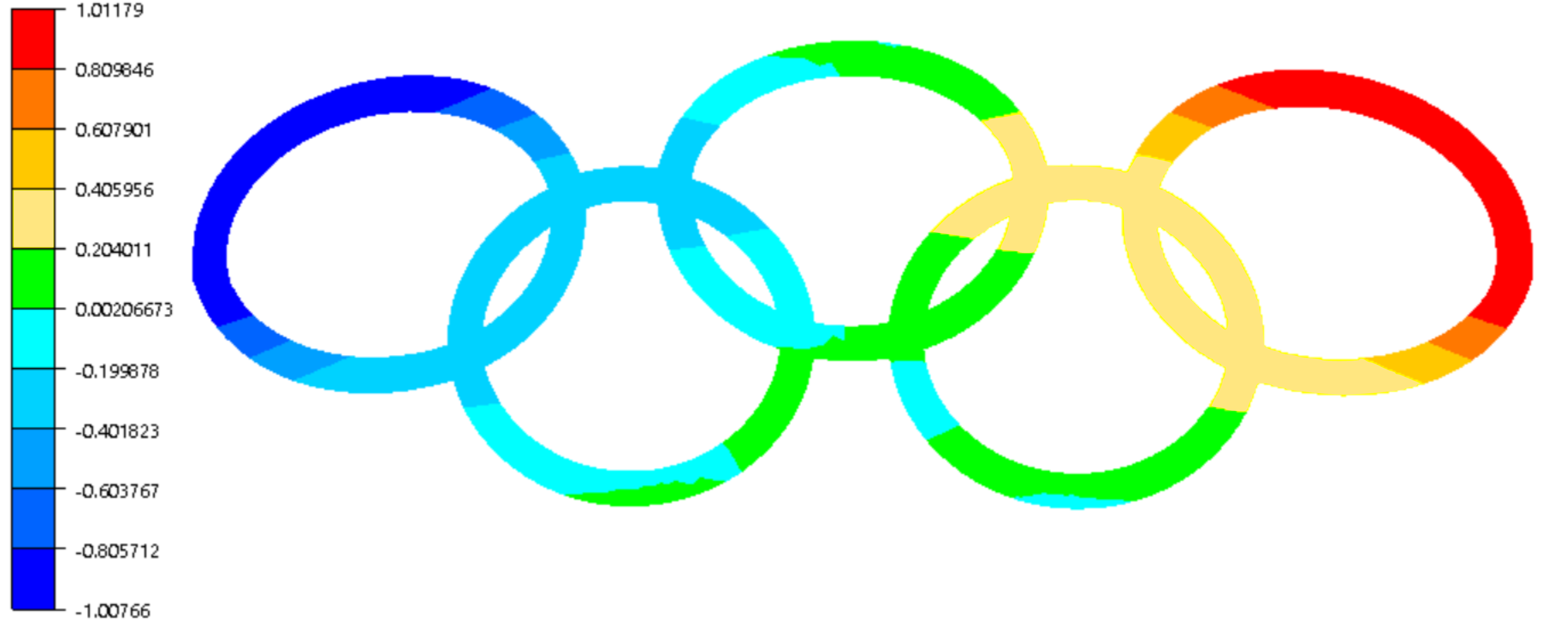}

        \end{subfigure}

        \begin{subfigure}[b]{0.55\linewidth}
            \caption{Predicted displacement}

            \includegraphics[width=\linewidth]{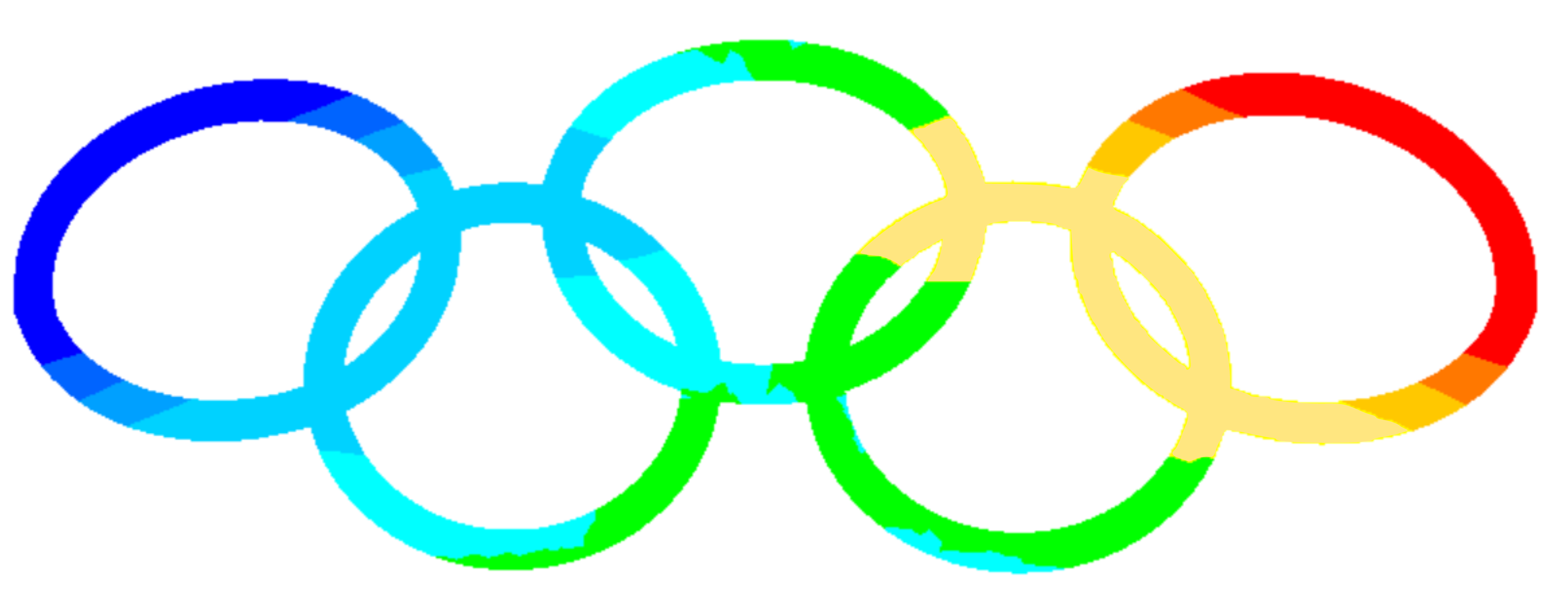}

        \end{subfigure}
    }

    \centerline{

        \rotatebox{90}{\hspace{1.2cm} $u_y$}

        \begin{subfigure}[b]{0.55\linewidth}
            \includegraphics[width=\linewidth]{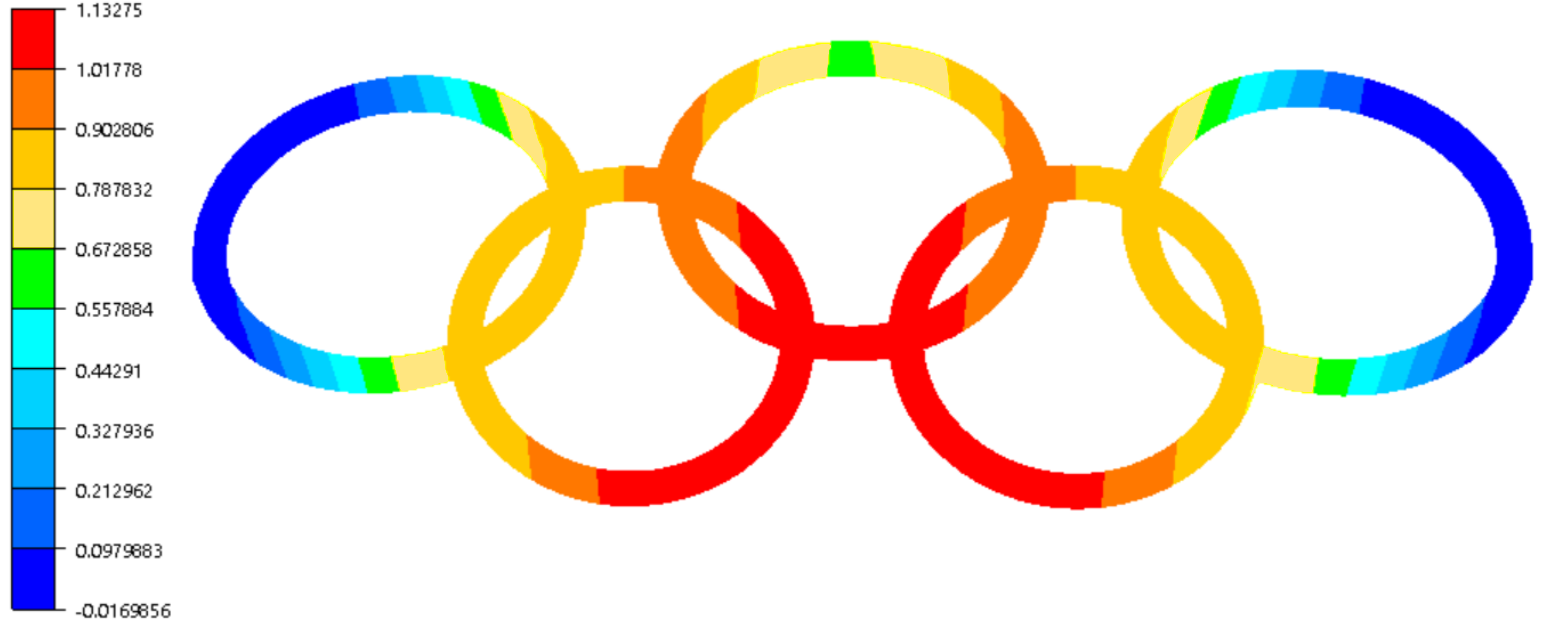}
        \end{subfigure}

        \begin{subfigure}[b]{0.55\linewidth}
            \includegraphics[width=\linewidth]{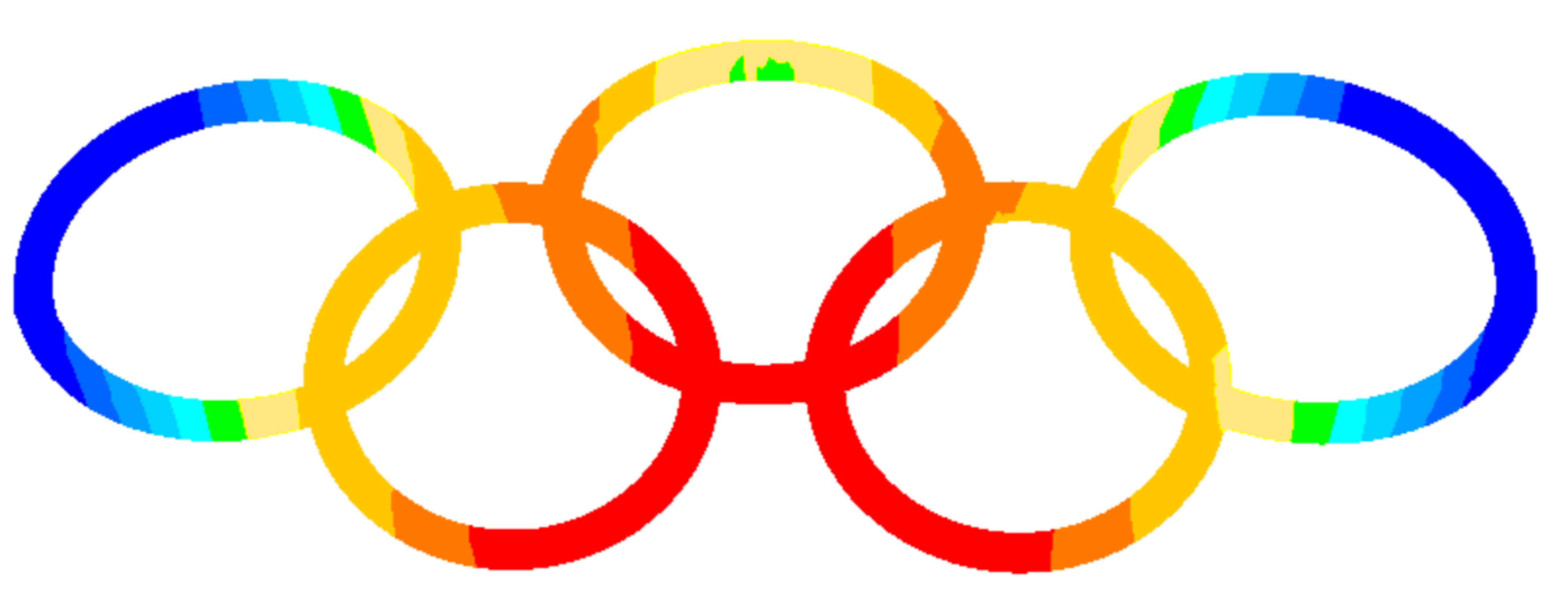}
        \end{subfigure}
    }
    \caption{Result of the training on the coarse mesh. (Left) Target displacement. (Right) Predicted displacement. (Top) $x$ component. (Bottom) $y$ component.}
    \label{fig:pred_true_coarse}\end{figure}

\begin{figure}
    \centerline{

        \rotatebox{90}{\hspace{1.2cm} $u_x$}

        \begin{subfigure}[b]{0.55\linewidth}
            \caption{Target displacement}
            \includegraphics[width=\linewidth]{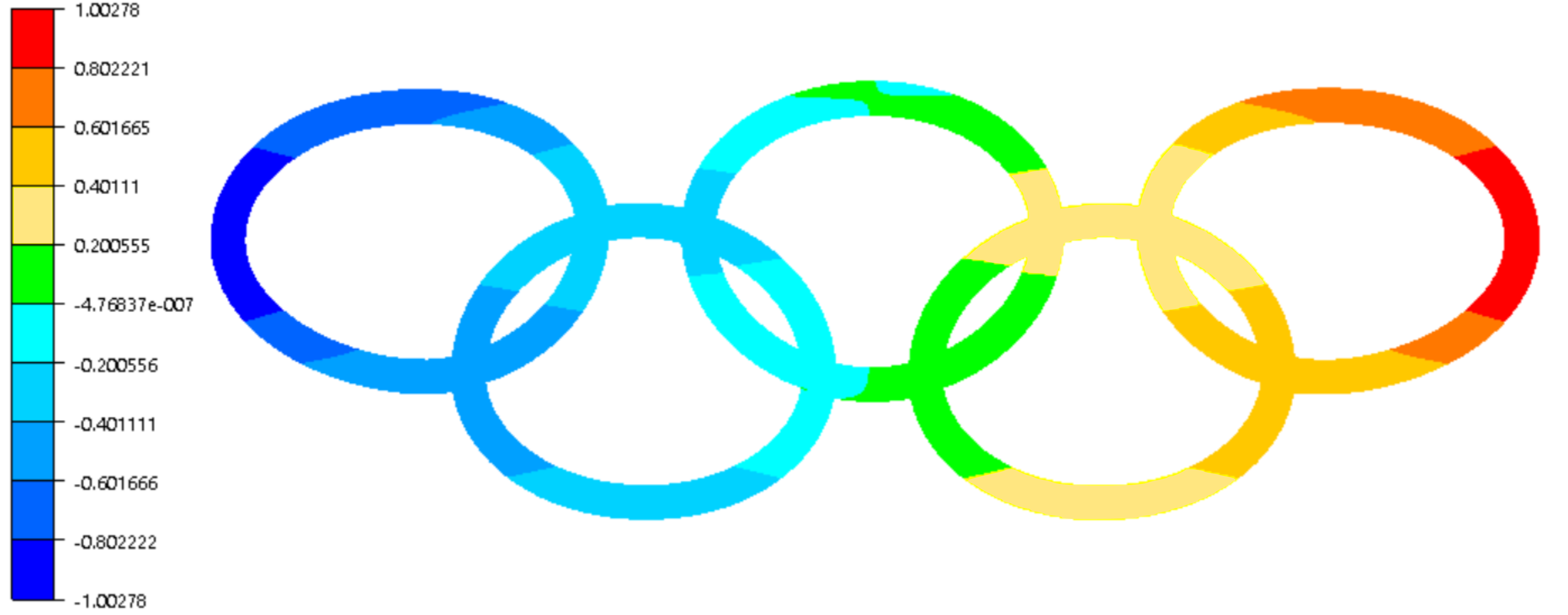}

        \end{subfigure}

        \begin{subfigure}[b]{0.55\linewidth}
            \caption{Predicted displacement}
            \includegraphics[width=\linewidth]{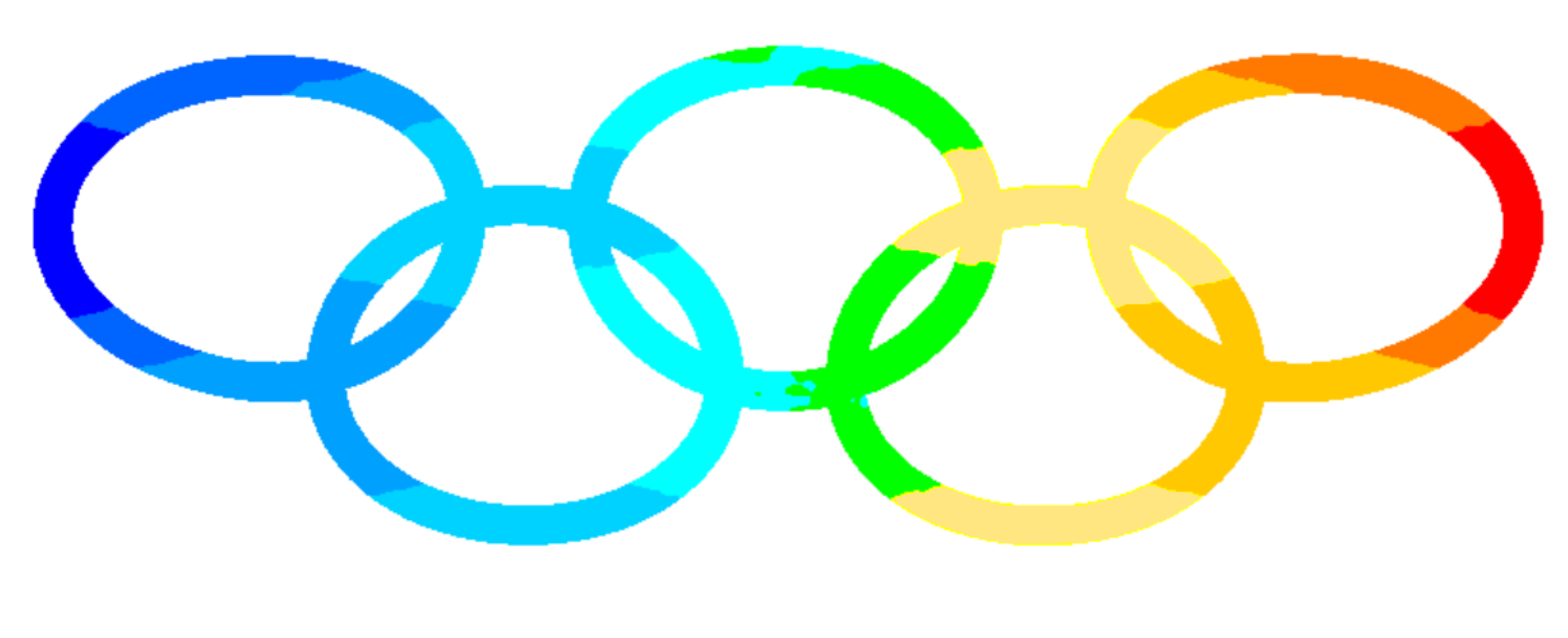}

        \end{subfigure}
    }

    \centerline{

        \rotatebox{90}{\hspace{1.2cm} $u_y$}

        \begin{subfigure}[b]{0.55\linewidth}
            \includegraphics[width=\linewidth]{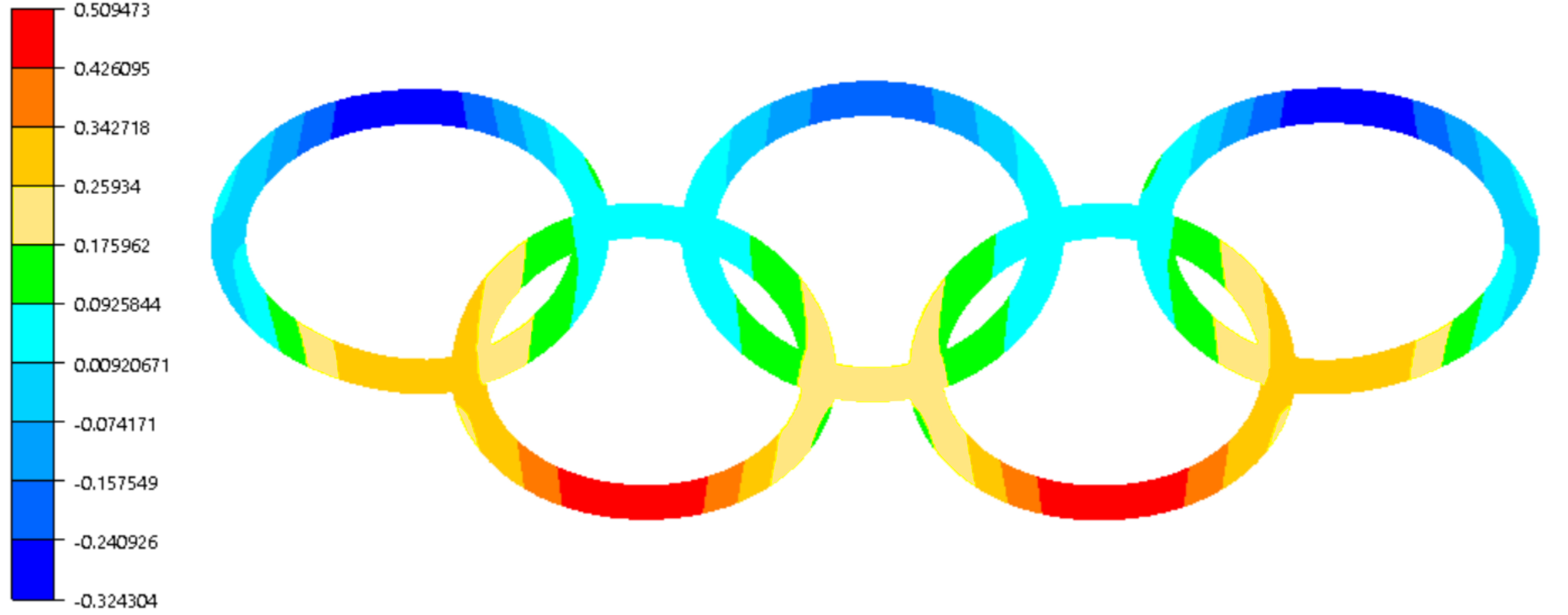}
        \end{subfigure}

        \begin{subfigure}[b]{0.55\linewidth}
            \includegraphics[width=\linewidth]{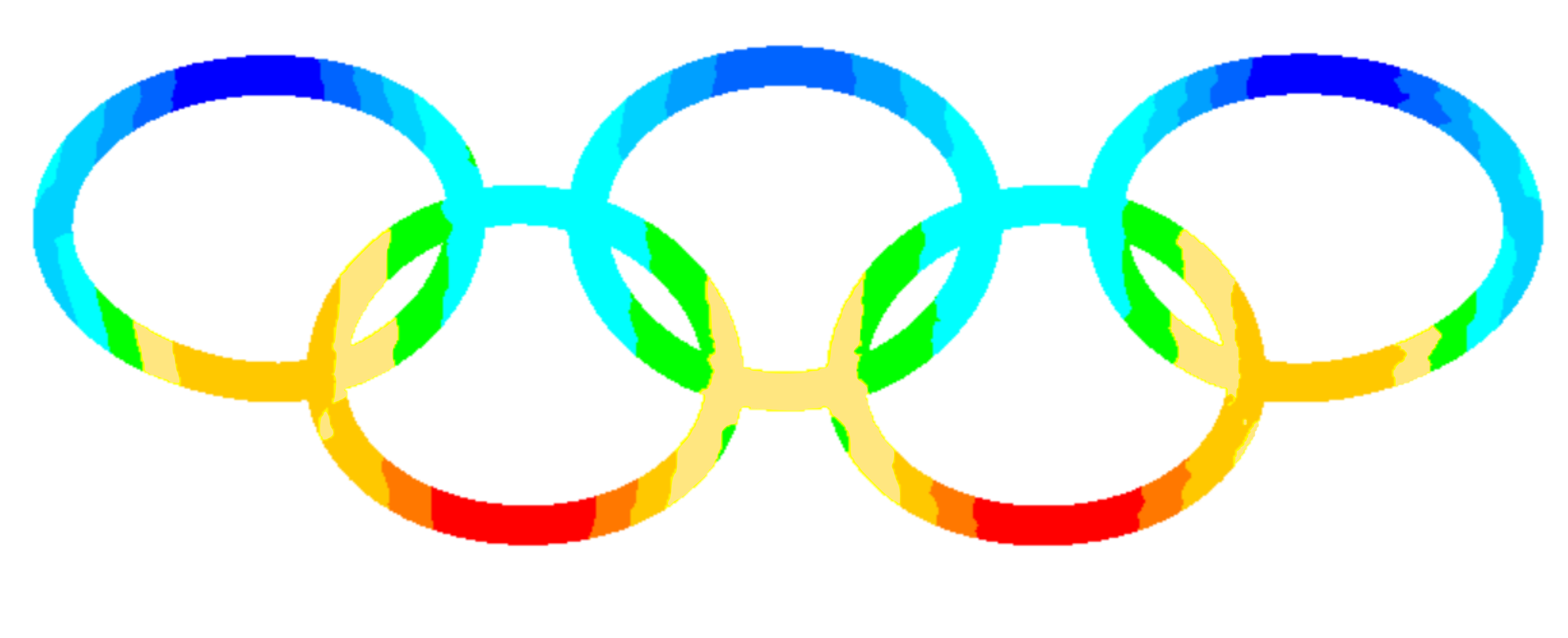}
        \end{subfigure}
    }
    \caption{Result of the training on the fine mesh. (Left) Target displacement. (Right) Predicted displacement. (Top) $x$ component. (Bottom) $y$ component.}
    \label{fig:pred_true_fine}\end{figure}

To perform the training, the same architecture as in section \ref{sec:model}, with adapted input size, has been used. As input to the model, the boundary conditions have been used (a null field, except on $\Gamma_1$ and $\Gamma_2$). The training has been performed with the Adam optimizer for 5000 epochs, with a learning rate of $5 \times 10^{-3}$. The training times are respectively $4.08 \times 10^2$ and $2.79 \times 10^3$ seconds on the coarse and fine mesh, on a single Intel Xeon Gold CPU. While combined optimizers, or a more thorough investigation of the training and model's hyperparameters, could yield better results, our main goal was to demonstrate the capacity of our model on various, complex cases, which are not adressed in the litterature. Therefore, an in-depth study of these hyperparameters has not been conducted on this case. Figures \ref{fig:pred_true_coarse} and \ref{fig:pred_true_fine} display the results of the training, respectively on the coarse and fine meshes.

Both figures show a good agreement between the predicted and true displacements, for the coarse and fine meshes. On the two cases, the relative mean squared errors between the predicted and true displacement fields are similar (approximately $1.1 \times 10^{-4}$). This result shows the good interpolation abilities of the physics-informed models. Unlike the traditional numerical approaches, coarse meshes may be sufficient to train accurate models, making this approach of interest. For an analysis of the extrapolation ability of the models, see section \ref{sec:numeric}.

\section{Numerical experiments in 3D}\label{sec:numeric}

\subsection{Presentation of the case}
\label{subsec:pb}

\textbf{PDE considered}

\noindent
To investigate the performance of the proposed model, a numerical study has been conducted. A three-dimensional electrostatic problem has been chosen. The unknown is the electric potential $V$ on a given, three-dimensional, irregular geometry $\Omega$. Two portions of the boundary, $\Gamma_1$ and $\Gamma_2$, are fixed, as boundary conditions, and the PDE that needs to be solved for $V$ is simply:

\begin{subequations}
    \begin{equation}
        \Delta V(x) = 0 \quad \forall \, x \in \Omega,
    \end{equation}

    \begin{equation}
        V(x) = 0 \quad \forall \, x \in \Gamma_1,
    \end{equation}

    \begin{equation}
        V(x) = 0.01 \quad \forall \, x \in \Gamma_2,
    \end{equation}
    \begin{equation}
        \frac{\partial V}{\partial n}(x) = 0 \quad \forall \, x \in \partial \Omega \backslash (\Gamma_1 \cup \Gamma_2).
    \end{equation}
    \label{eq:staticPb}
\end{subequations}

\noindent This problem has been solved using the commercial finite element solver FORGE \textregistered \cite{alves2017numerical}. In order to combine the classical FE solver with the physics-informed workflow, two key steps are distinguished. First, mesh and field data created by the software are used as input to construct the graph-based model. Then, the gradient kernel of the FE solver is used throughout the learning process, as detailed in section \ref{sec:customGrad}.  \newline

\noindent
\textbf{Geometry considered}

\begin{figure}[ht]
    \centering
    \includegraphics[width = 0.8\linewidth]{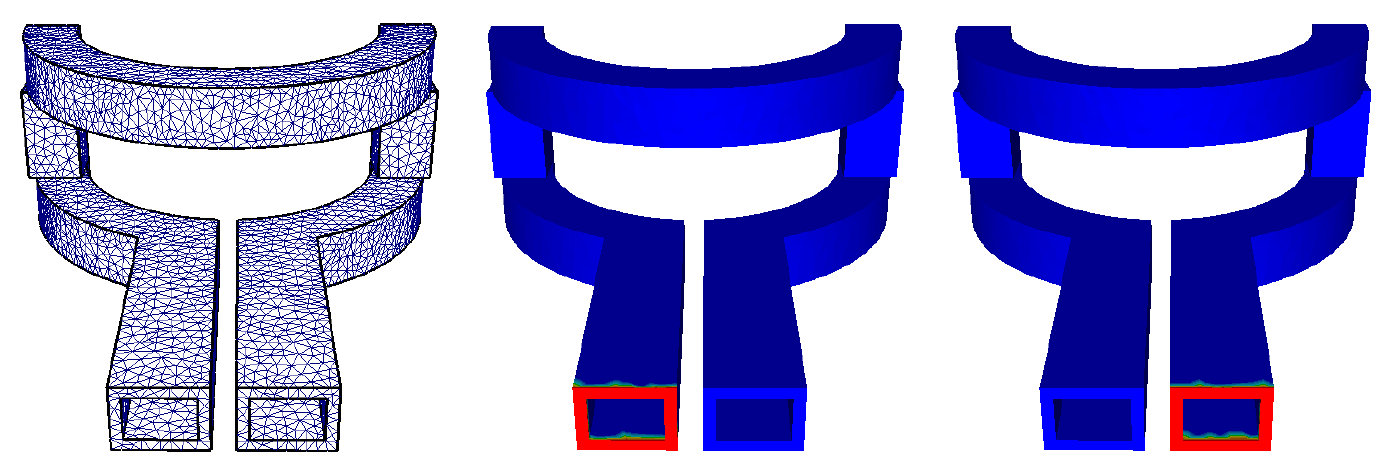}

    \caption{(Left) Geometry used for training, meshed domain. (Middle) Boundary set $\Gamma_1$, in red. (Right) Boundary set $\Gamma_2$, in red.}
    \label{fig:train_geometry}
\end{figure}

\begin{figure}[ht]
    \centering
    \includegraphics[width = \linewidth]{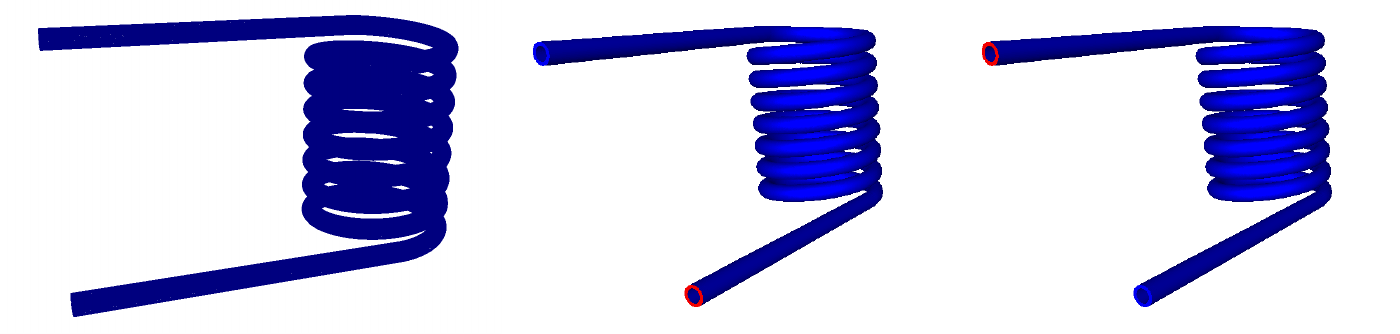}

    \caption{(Left) Geometry used for extrapolation, meshed domain. (Middle) Boundary set $\Gamma_1$, in red. (Right) Boundary set $\Gamma_2$, in red.}
    \label{fig:test_geometry}
\end{figure}

\noindent
One of the key aspects we want to investigate is the generalization capacity of the proposed model. Therefore, two different domains will be considered. The first one will be used for training, while the second one will be kept for the test step. For training, a domain with sharp angles has been selected. Its size is relatively small (10368 nodes) compared to the test domain size (485441 nodes). This discrepancy has been chosen to investigate the generalization capacity of the model, with respect to the domain size, but also its shape: the test domain is smoother, but also more complex, with many turns. The two meshed domains, along with the boundary sets, can be vizualized in figures \ref{fig:train_geometry} and \ref{fig:test_geometry}.

\begin{figure}[ht]
    \centering
    \includegraphics[width = 0.8\linewidth]{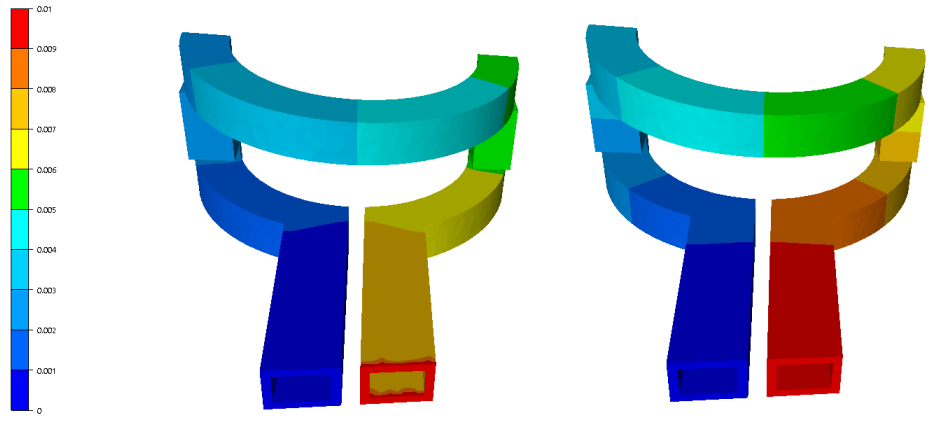}

    \caption{(Left) Approximation of the noisy FEM solution on a coarse mesh (used as input field to the model for training). (Right) FEM solution on the fine mesh, target field.}
    \label{fig:traincase}
\end{figure}

\begin{figure}[ht]
    \centering
    \includegraphics[width = 0.8\linewidth]{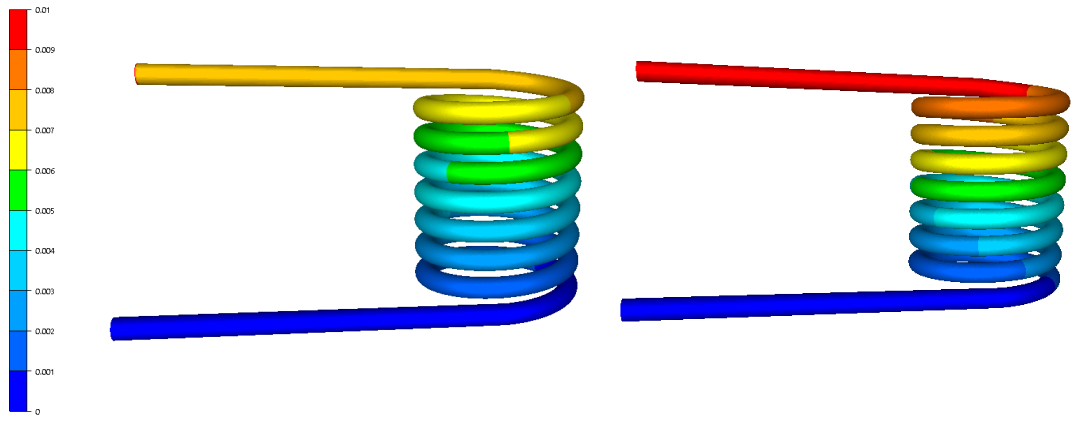}

    \caption{(Left) Approximation of the noisy FEM solution on a coarse mesh (used as input field to the model for extrapolation). (Right) FEM solution on the fine mesh, target field.}
    \label{fig:testcase}
\end{figure}

As an input to the model, for both geometries, we construct a noisy field with an approximate error of 20 \% compared to the true target field. This error is designed to simulate an approximate solution that could come from a very coarse finite element simulation for instance. Figures \ref{fig:traincase} and \ref{fig:testcase} present both cases considered, with the true fields and the input fields used during the train and test steps.

\subsection{Results on inference and generalisation}

\begin{figure}[ht]
    \centering
    \includegraphics[width = 0.8\linewidth]{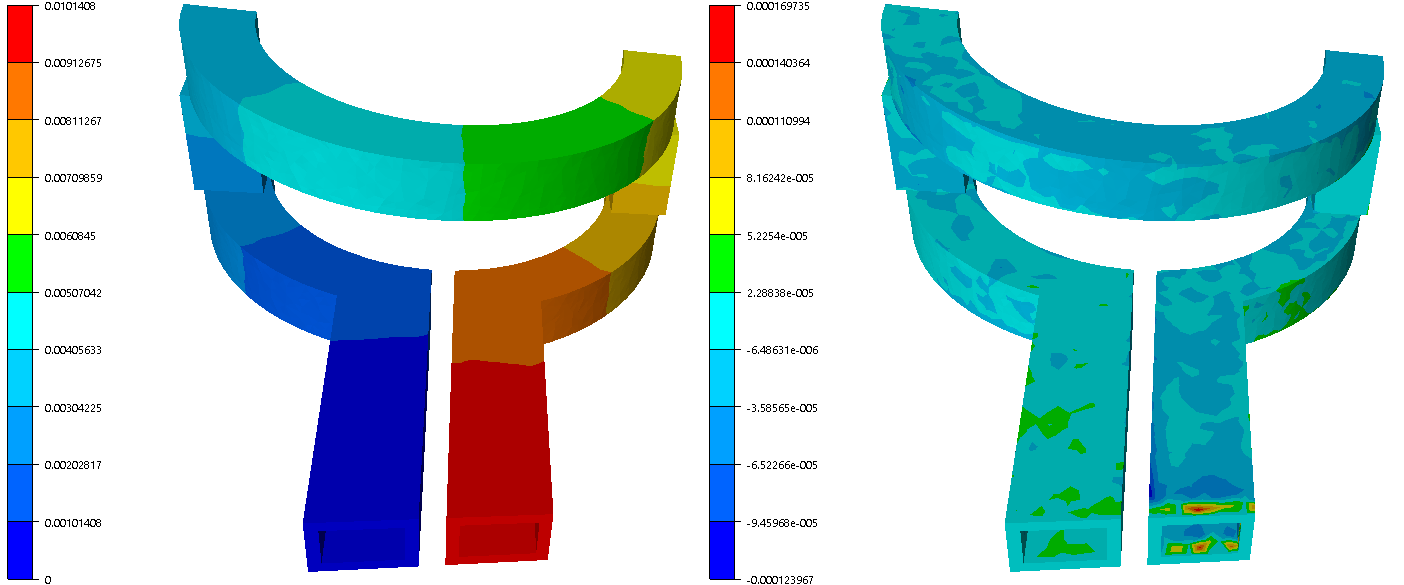}

    \caption{Result on the training case. Left: predicted potential. Right: error associated to the prediction. The relative $L^1$ error is $0.2 \%$.}
    \label{fig:predtrain}
\end{figure}

\begin{figure}[ht]
    \centering
    \includegraphics[width = 0.8\linewidth]{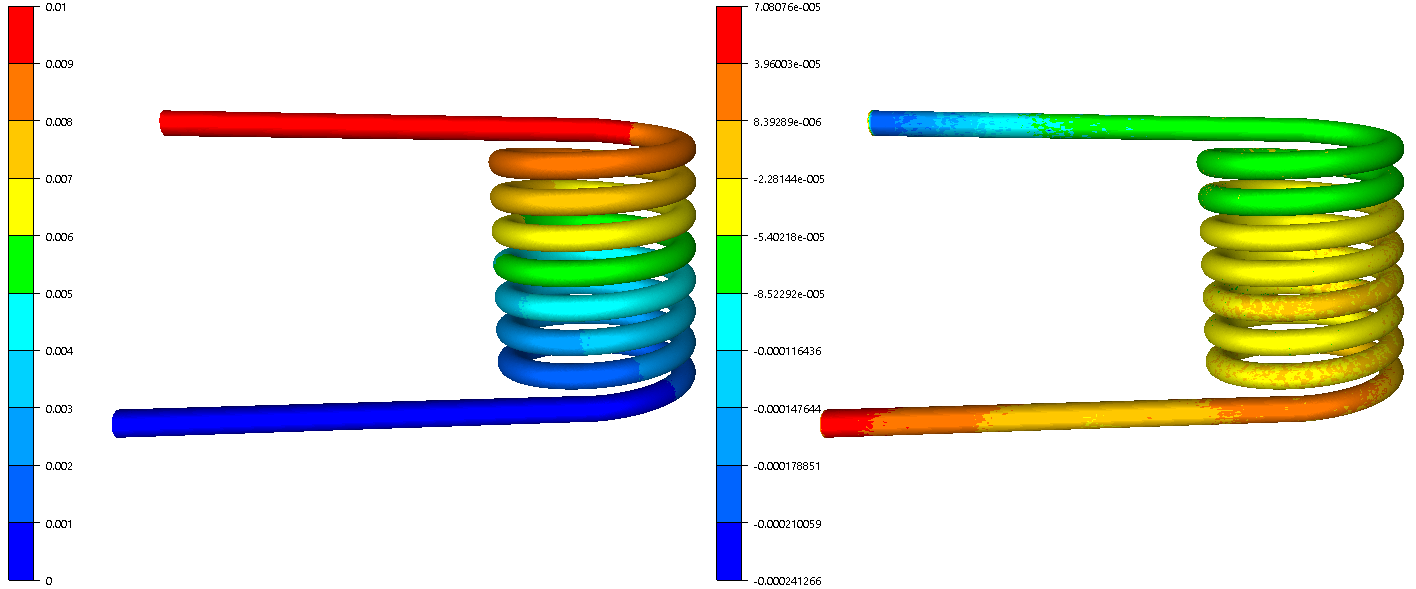}

    \caption{Result on the extrapolation case, with the model trained previously on the case described in Figure \ref{fig:traincase}. Left: predicted potential. Right: error associated to the prediction. The relative $L^1$ error is $0.7 \%$.}
    \label{fig:predtest}
\end{figure}

With the model architecture and the training procedure described before, a training has been conducted on the training geometry. For the training, a combination of 50 epochs of the Adam optimizer and 100 epochs of the L-BFGS optimizer has been selected.
Figures \ref{fig:predtrain} and \ref{fig:predtest} show the prediction of the model, along with the corresponding error, respectively for the training and the test case. The relative error is $0.2 \%$ for the training case, and $0.7 \%$ for the test case.

The results show a strong agreement between the prediction and the target, physical field, both in terms of accuracy and visually. Moreover, the results for the test case, although slightly worse than the training geometry, are still accurate and visually adequate. A pattern can be observed in the test geometry error, but the magnitude of the error is still significantly smaller than the magnitude of the true solution.

\subsection{Discussion on the model and the training procedure}

\subsubsection*{A comparison of the optimizers used on the training procedure}

\begin{figure}[ht]
    \centering
    \includegraphics[scale = 0.9]{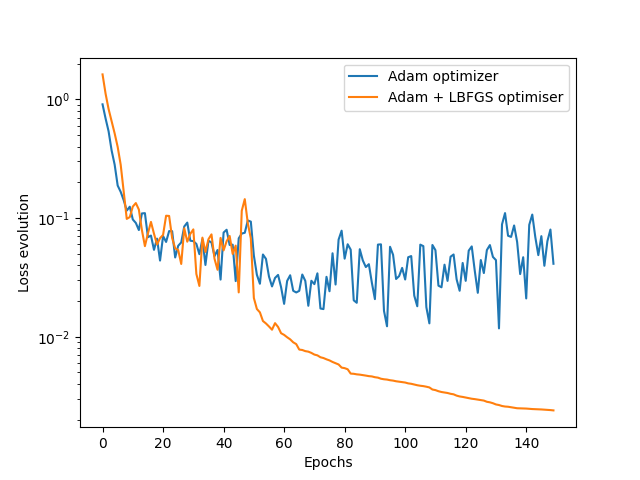}

    \caption{Evolution of the loss during training, for the Adam optimizer and for a combination of Adam and L-BFGS. For the Adam only training, the number of epochs is fixed to 150. For the combined training, 50 epochs are performed with Adam, and 100 with L-BFGS.}
    \label{fig:loss}
\end{figure}

As discussed, the best performance has been obtained with a combination of Adam and L-BFGS optimizers during training. Here, we plot the evolution of the loss during training with two strategies: a combination of Adam and L-BFGS, and using only Adam. Note that using only L-BFGS from the beginning led to a loss divergence. The result is presented in figure \ref{fig:loss}. The learning rate is fixed to $1\times10^{-4}$ for the Adam optimizer, and to $2 \times10^{-3}$ for L-BFGS.

The result clearly shows an improvement when using both optimizers combined, compared to performing only Adam. In both cases, the number of epochs has been fixed to 150 for a fair comparison.

\subsubsection*{Justification of the model: Ablation study}

We now present an ablation study to justify the impact of model's architecture and hybrid loss choices. The problem we adress is the one described in section \ref{subsec:pb}, with the training conducted on one domain, and the inference performed on the second domain. Both the training and test accuracies are reported. Two key aspects of our model have been investigated: the graph-based architecture and the loss formulation. Therefore, two other models have been selected. First, we train a model with the same architecture, but with a full data-driven approach, to demonstrate the performance of the chosen loss. This model is refered to as 'GraphNet, MAE train'. Then, we use our hybrid loss for a simple neural network, to demonstrate the performance of the chosen architecture. This model is refered to as 'NeuralNet'. The training was conducted for all the models with the same combination of Adam and L-BFGS. We found that this choice yielded the best results for the three models, even if in the case of the 'NeuralNet' model, there was no major difference between this training and an Adam-only optimizer. \\
'NeuralNet' has been selected to be a feedforward neural network, with two hidden layers and 128 neurons per layer. We found that adding hidden layers or increasing their width did not improve the accuracy of the model, which we explain by the fact that feedforward neural networks have poor inductive biases for this range of problems, compared to graph-based approaches. The results are presented in table \ref{table:ablation}.

\begin{table}[ht]
    \centering
    \begin{tabular}{|c||c|c|}
        \hline
        Model               & Relative train err. & Relative test err. \\
        \hline
        PiGMeN (ours)       & 0.2 \%              & \textbf{0.7 \%}    \\
        GraphNet, MAE train & \textbf{0.16 \%}    & 1.4 \%             \\
        NeuralNet           & 3.2 \%              & 11.9 \%            \\
        \hline
    \end{tabular}
    \caption{$L^1$ train and test errors for different versions of our model. PiGMeN is our full model. 'GraphNet, MAE train' refers to the same model's architecture, but trained only with the MAE loss. 'NeuralNet' refers to a vanilla neural network, trained with the hybrid loss described in section \ref{subsec:loss}.}
    \label{table:ablation}
\end{table}

This ablation study clearly shows the advantage of graph-based methods when dealing with complex domains, where the objective function has no simple algebraic expression. Moreover, direct coordinates are not suitable, since the problem is invariant in space: the feedforward neural network has poor generalization accuracy. Moreover, while further investigations should be conducted regarding the selected loss function, it seems that a physics-informed loss combined with a data-driven approach allows to better learn the underlying physical operator, which would explain the better generalization accuracy of our model, compared to a graph neural network trained on data only.

\section{Concluding remarks}

The main contributions of our work are threefold: first we have proven a limitation with the way PDE residuals are usually computed for PINNs, and we provided a way to overcome this limitation. Next, we proposed a novel model architecture, that respected the physical invariances of the problem. Finally, we proposed a hybrid loss, combining a data-driven approach and the physical knowledge we have on the case.

This work demonstrates the capability of graph neural networks trained with a physics-informed loss to handle complex, three-dimensional geometries, and to generalize to unseen settings. This performance is mainly due to overcoming a flaw in the standard method of computing PDE residuals in the physics-informed framework, and using an appropriate model architecture. The Lagrangian formulation was used to combine the physics-informed loss with available data. Numerical results demonstrate that our model, PiGMeN, is able to accurately solve an electrostatic problem on the training geometry, but also to generalize to completely different domains. To the best of our knowledge, such results had not been obtained before. \newline

For future works, further tests should be conducted to validate the approach, for instance with a full physics-informed approach, with no prior training data. To this end, an effort should be made to fully benefit from the Lagrangian formulation. An extension to time-dependent problems should also be considered. The proposed framework, which allows for calls to finite element solvers inside the deep learning loop, should enable significant progress in this direction.

\bibliographystyle{abbrv}
\bibliography{paper.bib}

\end{document}